# Towards Optimal Resource Allocation in Wireless Powered Communication Networks with Non-Orthogonal Multiple Access

Mariam M.N. Aboelwafa*, Mohamed A. Abd-Elmagid**, Alessandro Biason†, Karim G. Seddik*, Tamer ElBatt[II][‡] and Michele Zorzi†

*ECNG Dept., American University in Cairo, AUC Avenue, New Cairo 11835, Egypt
**Dept. of ECE, Virginia Tech, Blacksburg, VA 24060, USA
†DEI, University of Padova - via Gradenigo 6b, 35131 Padova, Italy
[II]Dept. of CSE, American University in Cairo, AUC Avenue, New Cairo 11835, Egypt
‡Dept. of EECE, Faculty of Engineering, Cairo University, Giza, Egypt
email: mariam.aboelwafa@aucegypt.edu, maelaziz@vt.edu, biasonal@dei.unipd.it,
kseddik@aucegypt.edu, telbatt@ieee.org, zorzi@dei.unipd.it

**Abstract**

The optimal allocation of time and energy resources is characterized in a Wireless Powered Communication Network (WPCN) with Non-Orthogonal Multiple Access (NOMA). We consider two different formulations; in the first one (max-sum), the sum-throughput of all users is maximized. In the second one (max-min), and targeting fairness among users, we consider maximizing the min-throughput of all users. Under the above two formulations, two NOMA decoding schemes are studied, namely, low complexity decoding (LCD) and successive interference cancellation decoding (SICD). Due to the non-convexity of three of the studied optimization problems, we consider an approximation approach, in which the non-convex optimization problem is approximated by a convex optimization problem, which satisfies all the constraints of the original problem. The approximated convex optimization problem can then be solved iteratively. The results show a trade-off between maximizing the sum throughout and achieving fairness through maximizing the minimum throughput.

## 1. Introduction

There has been a growing interest, recently, in studying new technologies for prolonging the lifetime of mobile devices [1]. RF Energy Harvesting (EH) is considered as a promising solution towards an unlimited power supply for wireless networks. However, it adds more complexity to system design and optimization [2], [3]. The main goal of this paper is to study an energy harvesting wireless network and characterize the optimal strategies to maximize the network throughput and ensure fairness among nodes.

There are two main paradigms in RF EH [4]; Simultaneous Wireless Information and Power Transfer (SWIPT) and Wireless Powered Communication Networks (WPCN).

In SWIPT, Wireless Energy Transfer (WET) and Wireless Information Transmission (WIT) occur simultaneously, in which energy and information are transmitted in the same signal [5]. In [6], Boshkovska designs a resource allocation algorithm for SWIPT systems. The algorithm design is formulated as a non-convex optimization problem for the maximization of the total harvested power at the EH receivers subject to quality of service (QoS) constraints. In [7], Ng and Schober study a resource allocation algorithm design for secure information and energy transfer to mobile receivers. In [8], multiple source-destination pairs communicate through their dedicated energy harvesting relays. A power splitting framework using game theory was developed to derive a profile of relays' power splitting ratios. Additionally, to overcome the problem that energy harvesting circuits are unable to harvest energy and decode information simultaneously, there are two proposed receiver designs in [9]: time switching and power splitting. By using the time switching setting, the receiving antenna periodically switches between energy harvesting and information decoding phases. On the other hand, under the power splitting, the received signal is split into two streams; one for the energy harvesting circuitry and the other is for information decoding. The application of SWIPT to Non-Orthogonal Multiple Access (NOMA) networks is investigated in [10], where Liu, et al. propose a new cooperative SWIPT-NOMA protocol in which users close to the source act as relays for far users' transmission.

In WPCN, the users harvest wireless energy from a dedicated Energy Rich (ER) source in the Downlink (DL), and then use it in the Uplink (UL) to send data to the Access Point (AP) [11]. In [12], it is assumed that the ER and the AP coincide. The optimal time allocations were characterized to achieve the max-sum throughput and the max-min throughput. A cooperative technique was studied in [13] and [14] to overcome the doubly near-far phenomenon. A WPCN with heterogeneous nodes (nodes with and without energy harvesting capabilities) was studied in [15] and it was shown how the presence of non-harvesting nodes can enhance the sum-throughput. [16] departed from the strong assumption adopted in [11] - [15], where the energy harvested in a slot is used completely in that slot, and, hence, embraces a long-term optimization framework. Additionally, [17] extended the long-term maximization of the half-duplex case in [16] to the full-duplex scenario. Conventional TDMA wireless networks were generalized in [18] to a new type of wireless networks: generalized-WPCNs (g-WPCNs), where nodes are equipped with RF energy harvesting circuitries along with energy supplies. It was shown that both conventional TDMA wireless networks and WPCNs with only RF energy harvesting nodes constitute lower bounds on the performance of g-WPCNs in terms of the max-sum throughput and max-min throughput.

Non-Orthogonal Multiple Access (NOMA) exploits an approach of user multiplexing in power domain [19], [20]. NOMA was introduced in WPCNs in [21] to enhance the power-bandwidth efficiency. It was shown in [22] that NOMA improves the spectral efficiency relative to orthogonal multiple access schemes. In [21], optimizing the time allocations was the main concern to maximize the sum-throughput of the slot-oriented case (all the harvested energy in a slot is also consumed in the same slot). Hence, Diamantoulakis, et al. introduced a sub-optimal policy for time allocations. Yuan and Ding investigated, in [23], the application of NOMA for the uplink (UL) of WPCNs. They maximize the sum rate by jointly designing the time allocation, the downlink (DL) energy beamforming and the receiver beamforming. In [24], two NOMA-based decoding schemes were introduced to maximize the sum-throughput of the network. Due to the difficulty of solving the optimization problem, an approximate iterative approach was proposed to solve a sub-problem and reach a sub-optimal solution. Chingoska and Nikoloska tackled, in [25], the doubly near-far effect in WPCNs by setting the decoding order signals received at the base station (BS) to be the inverse of the distances between the users and the BS.

This work considers the ER and AP as separate entities to accommodate a more general setting, in contrast to a large body of the literature where they coincide. Unlike [21], we jointly optimize time and power allocations over a finite horizon of $T > 1$ slots. In our system model, simultaneous transmissions can be successful if the received Signal-to-Interference-and-Noise Ratio (SINR) is higher than a pre-specified threshold, thus no accurate synchronization is required. Unlike [24], instead of solving the problem given one of the optimization decision variables, namely the allocated time for wireless energy transfer, we solve the original problem iteratively to attain performance superior to [24].

The main contributions of this work can be summarized as follows:

- Two decoding schemes are studied, namely: Low Complexity Decoding and Successive Interference Cancellation Decoding. The two schemes aim at optimizing the performance of a WPCN with and without interference cancellation.
- Maximizing the sum-throughput of the network is studied (max-sum). Since LCD leads to a non-convex problem, an iterative approach is introduced to solve two sub-problem in an alternating manner. On the other hand, the convexity of the max-sum problem with SICD is established and the problem is characterized to find the optimal transmission durations and powers.
- The fairness aspect is also studied and the optimization problem to maximize the minimum throughput of the network (for LCD and SICD) is characterized (max-min). Again, the problem is shown to be non-convex and an iterative algorithm is introduced to find an approximate solution that is close to the global optimum.

The rest of the paper is organized as follows. Section 2 presents the system model. The max-sum problem formulation for both LCD and SICD is presented in Section 3. To address the fairness issues of max-sum, the max-min problem formulation is presented in Sections 4. The numerical results are shown and discussed for key insights and observations in Section 5. Finally, Section 6 concludes the paper.

## 2. System Model

Consider a WPCN composed of one AP, one ER node and $K$ users $U_i, i = 1, 2, ..., K$ (see Figure 1). All nodes in the network are equipped with single antennas and operate over the same frequency. Only the ER node is equipped with a constant energy supply. It broadcasts DL wireless energy to the $K$ users in the network. Users receive energy and use the accumulated energy to send UL data to the AP.

Slotted time is considered and the slot duration is assumed to be normalized to one. Each slot $t = 1, 2, ..., T$ is divided into two phases: $\tau_{0,t}$ during which the ER broadcasts wireless energy on the DL to recharge the batteries of the devices and $1 - \tau_{0,t}$ during which all users transmit data to the AP independently and simultaneously over the UL. Note that all radios are half-duplex.

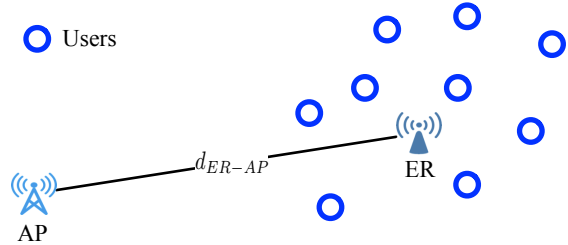

Figure 1: System Model

Users are uniformly distributed in a circle around ER at a distance $d_{U_i-ER}$, and $d_{U_i-AP}$ is the distance between $U_i$ and AP. Locations of all nodes are assumed to be known a priori, and therefore their average channel gains can be estimated. The DL channel power gain from ER to $U_i$ and the UL channel power gain from $U_i$ to AP, during time slot t, are denoted by $h_{i,t}$ and $g_{i,t}$, respectively. Hence, the harvested energy by $U_i$ in the DL phase can be expressed as [24]:

$$E_{i,t} = \eta_i h_{i,t} P_B \tau_{0,t} = \gamma_{i,t} \tau_{0,t}, \tag{1}$$

where $\eta_i$ denotes the energy harvesting circuitry efficiency [26], $P_B$ is the average transmit power by ER within $\tau_{0,t}$ and $\gamma_{i,t} \stackrel{\text{def}}{=} \eta_i h_{i,t} P_B$.

## 3. Max-Sum Throughput Optimization

In this section, the time and power allocation constrained optimization problem is formulated such that the sum UL throughput of the network is maximized. Two NOMA decoding schemes are studied, namely, Low Complexity Decoding (LCD) and Successive Interference Cancellation Decoding (SICD). The main objective is to formulate the problem of maximizing the achievable sum-throughput over a finite horizon of $T$ time slots subject to a number of constraints.

*3.1. Max-Sum Problem Formulation with Low Complexity Decoding*

Under the LCD scheme, the AP uses a single-user decoder to detect the signals received from all users without performing interference cancellation. Each signal suffers from interference from all other users. Thus, interference from other users is treated essentially as noise.

The achievable throughput of user $U_i$, in time slot $t$, can be expressed as:

$$R_{i,t} = (1 - \tau_{0,t}) \log_2(1 + x_{i,t}), \tag{2}$$

where $x_{i,t}$ is the average SINR at the AP for $U_i$ in time slot $t$.

In LCD, $x_{i,t}$ is given by:

$$x_{i,t} = \frac{g_{i,t}E_{i,t}}{\sigma^2(1-\tau_{0,t}) + \sum_{j=1,j\neq i}^{K} g_{j,t}E_{j,t}}, \quad (3)$$

where $E_{i,t}$ is the amount of energy used by $U_i$ in time slot $t$ and $\sigma^2$ is the noise power at the AP. It is worth noting that, under LCD, the interference term (i.e. summation) in the denominator includes interference from all users other than user $U_i$.

The sum throughput maximization problem can be formulated as:

$$\mathbf{P1_{LCD}}: \quad \max_{\tau_0, \mathbf{E}, \mathbf{x}} \sum_{t=1}^{T}\sum_{i=1}^{K} R_{i,t},$$

subject to: $Eq.(3)$,

$$\sum_{n=1}^{t} E_{i,n} \leq \sum_{n=1}^{t} \gamma_{i,n}\tau_{0,n}, \quad \forall t \quad \text{(energy causality constraints)},$$

$$x_{i,t} \geq S_i^{th}, \quad \forall i, \forall t \quad \text{(decoding constraints)},$$

$$0 \leq \tau_{0,t} \leq 1, \quad \forall t,$$

$$E_{i,t} \geq 0 \quad \forall i, \forall t,$$

where $\tau_0$, $\mathbf{E}$ and $\mathbf{x}$ are vectors whose elements are the harvesting time duration, the consumed energy by each user and the average SINR at the AP for each user over the finite horizon of $T$ slots, respectively. The role of the energy causality constraints are to guarantee that, in slot t, only the energy harvested in slots $\leq t$ can be used. The decoding constraints highlight the fact that if the SINR $x_{i,t}$ falls under a predefined threshold $S_i^{th}$, decoding will not be possible.

The objective function of $\mathbf{P1_{LCD}}$ is non-convex. Hence, the problem is a non-convex optimization problem. However, we exploit the problem structure and propose an efficient approach for solving this rather complex problem via splitting it into two separate subproblems, $\mathbf{P1_{LCD}}(\tau_0)$ and $\mathbf{P1_{LCD}}(E)$, that can be solved iteratively. The first is for a given harvesting slot duration $\tau_0$ and the later is for a given energy allocation vector $\mathbf{E}$. The proposed solution in this section is to solve these two sub-problems iteratively, alternating between a given ($\tau_0$) to get optimum energy allocation vector for the first sub-problem and then, for the obtained energy allocation vector, get the optimum harvesting slot duration in the second sub-problem, and so on alternating back and forth between the two sub-problems until convergence is attained.

*3.1.1. The First Sub-problem for a Given $\tau_0$: $\mathbf{P1_{LCD}}(\tau_0)$*

The first sub-problem given $\tau_0$ is formulated as:

$$\mathbf{P1_{LCD}}(\tau_0): \quad \min_{\mathbf{E},\mathbf{x}} \ \frac{1}{\prod_{t=1}^{T}\prod_{i=1}^{K}(1+x_{i,t})^{(1-\tau_{0,t})}},$$

$$\text{subject to:} \quad x_{i,t} \times \frac{\sigma^2(1-\tau_{0,t}) + \sum_{j=1,j\neq i}^{K} g_{j,t}E_{j,t}}{g_{i,t}E_{i,t}} \leq 1, \quad \forall i, \forall t,$$

$$\frac{\sum_{n=1}^{t} E_{i,n}}{\sum_{n=1}^{t} \gamma_{i,n}\tau_{0,n}} \leq 1, \quad \forall i, \forall t,$$

$$S_i^{th} x_{i,t}^{-1} \leq 1, \quad \forall i, \forall t,$$

$$E_{i,t} \geq 0 \quad \forall i, \forall t.$$

All constraints are expressed in the standard form for a Geometric Program (GP). However, the objective function of $\mathbf{P1_{LCD}}(\tau_0)$ is a ratio between two posynomials, thus, $\mathbf{P1_{LCD}}(\tau_0)$ is a nonconvex complementary GP [27]. Solving complementary GPs directly is NP-hard. Therefore, an approximate approach is used [24]. The denominator of the objective function, denoted $f(x)$, is approximated with a monomial function $\tilde{f}(x)$. In this case, the new approximate optimization problem becomes a standard GP that can be solved iteratively using standard techniques [28].

$\tilde{f}(x)$ is chosen to be:

$$\tilde{f}(x) = c \prod_{t=1}^{T} \prod_{i=1}^{K} (x_{i,t})^{y_{i,t}(1-\tau_{0,t})}, \tag{4}$$

where,

$$c = \frac{\prod_{t=1}^{T}\prod_{i=1}^{K}(1+\overline{x}_{i,t})^{y_{i,t}(1-\tau_{0,t})}}{\prod_{t=1}^{T}\prod_{i=1}^{K}(\overline{x}_{i,t})^{y_{i,t}(1-\tau_{0,t})}}, \tag{5}$$

$\overline{x}$ is the solution of the approximate GP in the previous iteration, and

$$y_{i,t} = \frac{\overline{x}_{i,t}}{1+\overline{x}_{i,t}}, \quad \forall i, \forall t. \tag{6}$$

Starting with an initial $\overline{x}$, we can obtain $c$ and $y_{i,t}$ from (5) and (6), respectively. With these values, we solve the approximate geometric program. The obtained solution can be used to get new values of $c$ and $y_{i,t}$. The procedure is repeated until the sum-throughput converges to a pre-specified accuracy.

*3.1.2. The Second Sub-problem for a Given $\mathbf{E}$: $\mathbf{P1_{LCD}}(E)$*

The second sub-problem given $\mathbf{E}$ is formulated as:

$$\mathbf{P1_{LCD}}(E): \qquad \max_{\boldsymbol{\tau_0},\mathbf{x}} \sum_{t=1}^{T}\sum_{i=1}^{K}(1-\tau_{0,t})\log_2(1+x_{i,t}),$$

subject to: *Eq.(3)*,

$$\sum_{n=1}^{t} E_{i,n} \leq \sum_{n=1}^{t} \gamma_{i,n}\tau_{0,n}, \quad \forall t \quad \text{(energy causality constraints)},$$

$$x_{i,t} \geq S_i^{th}, \quad \forall i, \forall t \quad \text{(decoding constraint)},$$

$$0 \leq \tau_{0,t} \leq 1, \quad \forall t.$$

Note that, by using (3), $x_{i,t}$ can be substituted and removed from the problem. All constraints are affine and it can be verified that the objective function has a Hessian matrix that is **positive semidefinite** (proof is omitted due to space limitations). Therefore, the problem $\mathbf{P1_{LCD}}(E)$ is convex and, hence, it can be solved using standard convex optimization tools.

The algorithm to solve $\mathbf{P1_{LCD}}$ iteratively, using the two sub-problems discussed above, is given in Algorithm 1.

---
**Algorithm 1** Solving $\mathbf{P1_{LCD}}$
---
1: **repeat**
2:     **procedure** SOLVE P1($\tau_0$)
3:         Initialize $\bar{x}$
4:         Compute $c$ and $y_{i,t}$ using (5) and (6)
5:         **repeat**
6:             Solve the approximate P1($\tau_0$)
7:             Update $c$ and $y_{i,t}$ using (5) and (6)
8:         **until** Sum throughput converges
9:         Find sub-optimum $R_{sum}$ and $\mathbf{E}$.
10:     **end procedure**
11:     **procedure** SOLVE P1($E$)
12:         Use standard convex optimization tools.
13:         Find sub-optimum $R_{sum}$ and $\boldsymbol{\tau_0}$.
14:     **end procedure**
15: **until** The maximized sum converges to a pre-specified accuracy
---

*3.2. Max-Sum Problem Formulation with Successive Interference Cancellation Decoding*

The LCD scheme, despite its highly complex computations (since it yields a non-convex optimization problem), it has simple implementation because interference is simply treated as noise. However, the performance is modest due to ignoring the structure of interference. This leads to a lower sum-throughput because each user suffers from interference from all other users. This motivates us to look at more sophisticated interference cancellation techniques to enhance the sum-throughput performance. In this sub-section, a Successive Interference Cancellation Decoding (SICD) scheme is introduced whereby interference can be partially cancelled out. The associated sum-throughput optimization problem is also formulated.

The SINR of $U_i$ in time slot $t$, after interference cancellation, can be expressed as:

$$x_{i,t} = \frac{g_{i,t}E_{i,t}}{\sigma^2(1-\tau_{0,t}) + \sum_{j=i+1}^{K} g_{j,t}E_{j,t}}, \quad \forall i. \tag{7}$$

It is worth noting that the SICD hinges on the assumption that every decoded single user's signal is removed from the interference term of all *next* users' signals to be decoded. Therefore, in (7), the interference term (i.e. summation) in the denominator includes interference from all users that will be decoded after the signal of user $U_i$. The achievable throughput for $U_i$ can be expressed as in (2) while substituting for the SINR $x_{i,t}$, from (7). Therefore, the achievable sum throughput over all $K$ users under SICD, which is independent of the users decoding order, can be expressed as (proof is omitted due to space limitations):

$$R_{sum}^{(t)} = (1-\tau_0)\log_2\left(1 + \frac{\sum_{j=1}^{K} g_i E_i}{\sigma^2(1-\tau_0)}\right). \tag{8}$$

The SICD sum-throughput optimization problem can then be formulated as:

$$\mathbf{P1_{SICD}}: \quad \max_{\tau_\mathbf{0},\mathbf{E},\mathbf{x}} \sum_{t=1}^{T}\sum_{i=1}^{K} R_{i,t},$$

$$\text{subject to:} \quad Eq.(7),$$

$$\sum_{n=1}^{t} E_{i,n} \leq \sum_{n=1}^{t} \gamma_{i,n}\tau_{0,n}, \quad \forall t \quad \text{(energy causality constraints)},$$

$$x_{i,t} \geq S_i^{th}, \quad \forall i, \forall t \quad \text{(decoding constraint)},$$

$$0 \leq \tau_{0,t} \leq 1, \quad \forall t,$$

$$E_{i,t} \geq 0 \quad \forall i, \forall t.$$

Recall that for a general function, $f(x)$, which is concave, its perspective function $g(x,t) = tf(\frac{x}{t})$ would also be concave [28]. By using $t = 1 - \tau_0^{(t)}$, the sum-throughput $R_{sum}^{(t)}$ is the perspective function of the concave function $\log_2\left(1 + \frac{\sum_{i=1}^{K} g_i^{(t)} E_i^{(t)}}{\sigma^2}\right)$. Therefore, $R_{sum}^{(t)}$ is a concave function in $[\tau_0^{(t)}, E_1^{(t)}, ..., E_K^{(t)}]$. Note that a non-negative weighted sum of concave functions is also concave, then the objective function of $\mathbf{P1_{SICD}}$ which is the non-negative weighted summation of $R_{sum}^{(t)}$, $\forall t$, is a concave function in $(\tau_0, \mathbf{E})$. In addition, all constraints of $\mathbf{P1_{SICD}}$ are affine in $(\tau_0, \mathbf{E})$. As a result, $\mathbf{P1_{SICD}}$ is a convex optimization problem, and, hence can be solved efficiently using standard convex optimization tools.

An efficient, yet simple, way to solve a constrained optimization problem, is to find its Lagrangian and solve the dual problem [28]. The Lagrangian of $\mathbf{P1_{SICD}}$ is given by:

$$\mathcal{L}(\mathbf{E}, \tau_0, \lambda, \mu) = \sum_{t=1}^{T}\sum_{i=1}^{K} R_{i,t} + \sum_{i=1}^{K}\sum_{n=1}^{T} \lambda_{i,n}\left(\sum_{t=1}^{n}\left(\gamma_{i,t}\tau_{0,t} - E_{i,t}\right)\right)$$
$$+ \sum_{i=1}^{K}\sum_{t=1}^{T} \mu_{i,t}\left(x_{i,t} - S_i^{th}\right),$$

where $\lambda_{i,t}$ and $\mu_{i,t}$ are the dual variables associated with the energy causality and practical decoding constraints. Now, we need to solve the following optimization problem (dual problem), namely $\mathbf{D}_{SICD}$:

$$\mathbf{D}_{SICD}: \quad \max_{\tau_0, \mathbf{E}} \quad \mathcal{L}(\mathbf{E}, \tau_0, \lambda, \mu),$$
$$\text{subject to:} \quad 0 \leq \tau_{0,t} \leq 1, \quad \forall t,$$
$$E_{i,t} \geq 0 \quad \forall i, \forall t.$$

*Lemma:* Given $\lambda$ and $\mu$, the optimal time and energy allocations of $\mathbf{D_{SICD}}$ are given by:

$$\tau_{0,t}^* = \min\left[\left(1 - \frac{\sum_{i=1}^{K} g_{i,t} E_{i,t}}{z_t^* \sigma^2}\right)^+, 1\right]. \tag{9}$$

$$E_{i,t}^* = \left(\frac{(1-\tau_{0,t})(g_{i,t} - \sigma^2 a_{i,t})}{a_{i,t} g_{i,t}} - \frac{1}{g_{i,t}}\sum_{j \geq i}^{K} g_{i,t} E_{i,t}\right)^+. \tag{10}$$

The variable $a_{i,t}$ is defined as:

$$a_{i,t} \stackrel{\text{def}}{=} \ln(2)\left(\sum_{n=t}^{T} \lambda_{i,n} + g_{i,t}\chi\{i \geq 2\}\sum_{j=1}^{i-1} \mu_{j,t} S_j^{th} - \mu_{i,t} g_{i,t}\right), \tag{11}$$

where $\chi\{.\}$ is the indicator function and $(.)^+ \stackrel{\text{def}}{=} \max\{0,.\}$ and $z_t^*$ is the unique solution of $f(z_t) = b^{(t)}$, where $f(z)$ and $b^{(t)}$ are given by:

$$f(z_t) = \ln(1 + z_t) - \frac{z_t}{1 + z_t}. \tag{12}$$

$$b^{(t)} = \ln(2)\left(\sigma^2 \sum_{i=1}^{K} \mu_{i,t} S_i^{th} + \sum_{i=1}^{K}\sum_{n=t}^{T} \lambda_{i,n} \gamma_{i,t}\right). \tag{13}$$

*Proof:* It can be verified that there exists $\boldsymbol{\tau_0}$ and $\mathbf{E}$ that strictly satisfy all the constraints of $\mathbf{D_{SICD}}$. Hence, strong duality holds for this problem [28]; therefore, the KKT conditions given below are necessary and sufficient for the global optimality:

$$\frac{\delta}{\delta \tau_{0,t}} \mathcal{L} = \ln\left(1 + \frac{\sum_{j=1}^{K} g_{j,t} E_{j,t}}{\sigma^2(1 - \tau_{0,t})}\right) - \frac{\sum_{j=1}^{K} g_{j,t} E_{j,t}}{\sigma^2(1 - \tau_{0,t}) + \sum_{j=1}^{K} g_{j,t} E_{j,t}} - a^{(t)} = 0. \quad (14)$$

$$\frac{\delta}{\delta E_{i,t}} \mathcal{L} = \frac{\frac{g_{i,t}}{\sigma^2}}{1 + \frac{\sum_{i=1}^{K} g_{i,t} E_{i,t}}{\sigma^2(1 - \tau_{0,t})}} - b_{i,t} = 0, \quad (15)$$

$\forall i$ and $t$, where $a_{i,t}$ and $b^{(t)}$ are given by (11) and (13), respectively.

By defining $z_t = \frac{\sum_{j=1}^{K} g_{i,t} E_{j,t}}{\sigma^2(1 - \tau_{0,t})}$, (14) can be reformulated as $f(z_t) = b^{(t)}$, where $f(z_t)$ is given in (12). Since $f(z_t)$ can be verified to be a monotonically increasing function of $z_t \geq 0$, where $f(0) = 0$, then there exists a unique solution $z_t^*$ that satisfies $f(z_t^*) = b^{(t)}$ and, hence, $\tau_{0,t}^*$ can be expressed as in (9) and by using (15), $E_{i,t}^*$ can be expressed as in (10).

## 4. Max-Min Throughput Optimization

In this section, we address the potential unfairness typically exhibited by sum-throughput optimal policies. Under the max-sum throughput formulation, some nodes are likely to be allocated very little, or no, resources (time and power) in some scenarios, such that they achieve almost zero throughput. This can be a serious problem for some applications. In wireless sensor networks, for instance, in which all sensors need to periodically send their sensing data to the AP with the same rate, fairness is a necessity. To address fairness in our problem context, we adopt a widely accepted sense of fairness in the communications and networks literature, that is, max-min fairness [29], [30]. Next, the problem formulation and solution approach for maximizing the minimum UL throughput in LCD and SICD is discussed in details.

*4.1. Max-Min Problem Formulation with Low Complexity Decoding*

The max-min UL throughput problem, with low complexity decoding, can be formulated as follows:

$$\mathbf{P2_{LCD}}: \quad \max_{\boldsymbol{\tau_0},\mathbf{E},\mathbf{x},\overline{\mathbf{R}}} \overline{R},$$

$$\text{subject to:} \quad R_{i,t} \geq \overline{R},$$
$$Eq.(3),$$
$$\sum_{n=1}^{t} E_{i,n} \leq \sum_{n=1}^{t} \gamma_{i,n}\tau_{0,n}, \quad \forall t \quad \text{(energy causality constraints)},$$
$$x_{i,t} \geq S_i^{th}, \quad \forall i, \forall t \quad \text{(decoding constraint)},$$
$$0 \leq \tau_{0,t} \leq 1, \quad \forall t,$$
$$E_{i,t} \geq 0 \quad \forall i, \forall t,$$

where $\overline{R}$ is the minimum throughput to be maximized, $R_{i,t}$ is the achievable UL throughput of $U_i$ in time slot $t$ and is expressed in (2) and $x_{i,t}$ is the average SINR at the AP for $U_i$ in time slot $t$.

This problem differs from the two max-sum problems, studied in Section 3, in the objective function and the first constraint. The objective function is affine but the first constraint is not convex. Hence, $\mathbf{P2_{LCD}}$ is a non-convex problem.

In order to circumvent the non-convexity hurdle of $\mathbf{P2_{LCD}}$, we adopt an iterative solution approach similar to the one followed in Section 3.1 to solve $\mathbf{P1_{LCD}}$. To this end, we split $\mathbf{P2_{LCD}}$ into two sub-problems $\mathbf{P2_{LCD}}(\boldsymbol{\tau_0})$ (given $\boldsymbol{\tau_0}$) and $\mathbf{P2_{LCD}}(\mathbf{E})$ (given $\mathbf{E}$). The first sub-problem is still non-convex and will be solved using an approximate iterative method. On the other hand, the second sub-problem is convex and can be solved using standard convex optimization tools.

*4.1.1. The First Sub-problem for a give $\boldsymbol{\tau_0}$: $\mathbf{P2_{LCD}}(\tau_0)$*

The first sub-problem given $\boldsymbol{\tau_0}$ is formulated as follows:

$$\max_{\mathbf{E},\mathbf{x},\overline{\mathbf{R}}} \overline{R},$$

$$\text{subject to:} \quad R_{i,t} \geq \overline{R} \quad \forall i, \forall t,$$
$$Eq.(3),$$
$$\sum_{n=1}^{t} E_{i,n} \leq \sum_{n=1}^{t} \gamma_{i,n}\tau_{0,n}, \quad \forall t \quad \text{(energy causality constraints)},$$
$$x_{i,t} \geq S_i^{th}, \quad \forall i, \forall t \quad \text{(decoding constraint)},$$
$$E_{i,t} \geq 0 \quad \forall i, \forall t.$$

The objective function and the constraints are affine, except for the first constraint which is non-convex. To solve this problem, an approximate iterative approach is adopted [31].

The approach used to solve $\mathbf{P2_{LCD}}(\tau_0)$ is based on solving a sequence of strongly convex inner approximations of the problem until a stationary solution of $\mathbf{P2_{LCD}}(\tau_0)$ is reached. This solution guarantees, based on the proof and assumptions in [31], the feasibility of the solutions in every iteration.

The approach in [31] is based on replacing a non-convex objective function (say $U(\mathbf{x})$) by a strongly convex and simple function ($\tilde{U}(\mathbf{x};\mathbf{y})$) and constraints ($g_m(\mathbf{x})$, where $m$ is the non-convex constraints index) with convex upper estimates ($\tilde{g}_m(\mathbf{x};\mathbf{y})$) to create a sub-problem $\mathbf{P_y}$. The sub-problem $\mathbf{P_y}$ is strongly convex and has a unique solution $\hat{\mathbf{x}}(\mathbf{y})$ (a function of $\mathbf{y}$). By starting from a feasible point $\mathbf{y}^{(0)}$, the proposed method, iteratively, computes the solution of the sub-problem $\mathbf{P_y}$, which is $\hat{\mathbf{x}}(\mathbf{y})$ and then takes a step ($\zeta^n \in (0,1]$, where $n$ is the iteration index) from $\mathbf{y}$ towards $\hat{\mathbf{x}}(\mathbf{y})$.

Note that the point $\mathbf{y}$ generated by the algorithm in every iteration is always feasible for the original problem $\mathbf{P2_{LCD}}(\tau_0)$. Convergence is guaranteed under mild assumptions that offer a lot of flexibility in the choice of the approximation functions and free parameters.

The main problem that affects this approach is the *affine* objective function. To check stationarity of every iteration, we need the objective function to be a function of $\mathbf{y}$ to study its gradient until a stationary solution is reached. To solve this problem, the objective function $\overline{R}$ is replaced by a single user throughput (the one that had the minimum value at the initial search point: $R_{l,t}$) and an equality constraint is added to attain the same target. The objective function must be modified every iteration based on the minimum achieved value at the new initials (the optimum point of the previous iteration).

After this modification, problem $\mathbf{P2_{LCD}}(\tau_0)$ can be formulated as follows:

$$\max_{\mathbf{E},\mathbf{x},\overline{\mathbf{R}}} \quad R_{l,t},$$

$$\begin{aligned}
\text{subject to:} \quad & R_{l,t} = \overline{R} \quad \forall t, \\
& R_{i,t} \geq \overline{R} \quad \forall i = 1,2,3,...,K, i \neq l, \forall t, \\
& Eq.(3), \\
& \sum_{n=1}^{t} E_{i,n} \leq \sum_{n=1}^{t} \gamma_{i,n} \tau_{0,n}, \quad \forall t \quad \text{(energy causality constraints)}, \\
& x_{i,t} \geq S_i^{th}, \quad \forall i, \forall t \quad \text{(decoding constraint)}, \\
& E_{i,t} \geq 0 \quad \forall i, \forall t.
\end{aligned}$$

The chosen approximation is given by (16) [31]:

$$\tilde{U}(x;y) = \nabla_x U(y)^T (x-y) + \frac{1}{2}|x-y|^2. \tag{16}$$

This approximation mimics proximal gradient methods. It can be used if no convexity whatsoever is present. It is also used for the first two constraints, but we must check that the approximation $\tilde{g}_m(\mathbf{x}; \mathbf{y})$ is an upper approximate function if the non-convex constraint is $g_m(\mathbf{x}) \leq C$, where $C$ is a constant and a lower approximate function if $g_m(\mathbf{x}) \geq C$.

*4.1.2. The Second Sub-problem for a given* $\mathbf{E}$: $\mathbf{P2_{LCD}}(E)$

The second sub-problem given $\mathbf{E}$ is formulated as follows:

$$\max_{\tau_0, \mathbf{x}, \overline{\mathbf{R}}} \overline{R},$$

$$\text{subject to:} \quad R_{i,t} \geq \overline{R},$$
$$Eq.(3),$$
$$\sum_{n=1}^{t} E_{i,n} \leq \sum_{n=1}^{t} \gamma_{i,n} \tau_{0,n}, \quad \forall t \quad \text{(energy causality constraints)},$$
$$x_{i,t} \geq S_i^{th}, \quad \forall i, \forall t \quad \text{(decoding constraint)},$$
$$0 \leq \tau_{0,t} \leq 1, \quad \forall t.$$

The objective function in addition to all constraints are affine except for the first constraint. It can be verified that the left hand side of the first constraint has a Hessian matrix that is **positive semidefinite**, hence, it is convex. Therefore, problem $\mathbf{P2_{LCD}}(\mathbf{E})$ is convex and can be solved efficiently using standard convex optimization tools.

*4.2. Max-Min Problem Formulation with Successive Interference Cancellation Decoding*

The max-min optimization problem with SICD is formulated in a similar manner as LCD except for minor differences:

$\mathbf{P2_{SICD}}$ : $\max_{\tau_0, \mathbf{E}, \mathbf{x}, \overline{\mathbf{R}}} \overline{R},$

$$\text{subject to:} \quad R_{i,t} \geq \overline{R},$$
$$Eq.(7),$$
$$\sum_{n=1}^{t} E_{i,n} \leq \sum_{n=1}^{t} \gamma_{i,n} \tau_{0,n}, \quad \forall t \quad \text{(energy causality constraints)},$$
$$x_{i,t} \geq S_i^{th}, \quad \forall i, \forall t \quad \text{(decoding constraint)},$$
$$0 \leq \tau_{0,t} \leq 1, \quad \forall t,$$
$$E_{i,t} \geq 0 \quad \forall i, \forall t.$$

Notice that the only difference between $\mathbf{P2_{LCD}}$ and $\mathbf{P2_{SICD}}$ is the definition of the SINR which is denoted by $x$. As mentioned earlier, in LCD, the interference term in (3) includes signals from all other users. On the other hand, in SICD, it includes only the signals which has not been decoded yet and therefore are not cancelled from interference, as given in (7).

Under the SICD scheme, it is important to notice that although the decoding order, that is, the order by which the users' signals are decoded at the AP, doesn't affect the sum throughput of the network as shown in [32], it will certainly affect the fairness aspect. To achieve the highest fairness, it can be easily proven that the optimal decoding order w.r.t fairness is based on the received signal strength that is directly affected by the channel gain; from the strongest to the weakest signal [25].

The approach used to solve **P2$_{\mathbf{SICD}}$** is the same used to solve **P2$_{\mathbf{LCD}}$**. It is worth mentioning that, for **P2$_{\mathbf{LCD}}$** and **P2$_{\mathbf{SICD}}$**, if the optimum solution can be reached, the max-min throughput would have been the common throughput achieved by all users. However, since an approximate approach is used, the users do not reach a common throughput and some variations can be noticed between the throughputs achieved by individual users.

The formal description of the approach to solve **P2$_{\mathbf{LCD}}$** and **P2$_{\mathbf{SICD}}$** is given in Algorithm 2.

---

**Algorithm 2** Solving **P2$_{\mathbf{LCD}}$** and **P2$_{\mathbf{SICD}}$**

---
1: **repeat**
2:    **procedure** SOLVE P2($\tau_0$)
3:       Initialize $n = 0, \zeta^{(n)} \in (0,1], y^{(n)} \in$ feasible set.
4:       **repeat**
5:          Choose the objective function to be the user throughput with minimum value at initial point.
6:          Approximate non-convex function and constraints.
7:          Compute $y^{(n)}$, the solution of the sub-problem $\mathbf{P_y}$.
8:          Set $y^{(n+1)} = y^{(n)} + \zeta^{(n)}(\hat{x}(y^{(n)}) - y^{(n)})$.
9:          $n \leftarrow n + 1$.
10:      **until** Stationary solution of P2($\tau_0$) is reached
11:      Find sub-optimum $R_{min}$ and **E**.
12:    **end procedure**
13:    **procedure** SOLVE P2($E$)
14:       Use standard convex optimization tools.
15:       Find sub-optimum $R_{min}$ and $\tau_0$.
16:    **end procedure**
17: **until** The maximized throughput converges to a pre-specified accuracy

---

## 5. Performance Evaluation

*5.1. Simulation Setup*

We consider $K$ single-antenna nodes, where $K$ ranges from 2 to 20. Nodes are distributed randomly around the ER node in a circular area of radius 10 meters. A horizon of $T$ timeslots, ranging from 1 to 10, is studied. Nodes receive DL power from the ER node and send UL data to the AP which is located $d_{ER-AP}$ meters from the ER. For the UL transmission, the path loss model is $g_{i,t} = 10^{-3} d_{U_i-AP}^{-2}$. For the DL, we use the parameters of the P2110 device [26]. The system parameters, used to generate numerical results, are listed in Table 1.

Table 1: Simulation Parameters

| Parameters | Definition | Values |
|---|---|---|
| $d_{ER-AP}$ | Distance between ER and AP (meters) | 0:20:120 |
| $\sigma^2$ | Noise Power | -155 dBm/Hz |
| $BW$ | Bandwidth | 1 MHz |
| $P_B$ | ER Node Transmission Power | 3 W |
| $f_c$ | Central Frequency | 915 MHz |
| $G_r$ | Receiver Antenna Gain | 6 dB |
| $\eta_i$ | Harvesting Efficiency | 0.49 |
| $S_i^{th}$ | Decoding Threshold (dB) | -10:1:0 |

*5.2. Numerical Results*

Simulations were carried out using the optimization toolbox in MATLAB. The performance results presented next revolve around three main thrusts, namely max-sum optimization, max-min optimization and the fundamental throughput-fairness trade-off within our problem context.

*5.2.1. Max-Sum Performance Results*

Recall that, under the max-sum problem formulation, the objective function to be maximized is the sum UL throughput. The performance of the presented approaches, namely LCD and SICD, are presented next.

We compare three solutions for the max-sum problem under LCD (**P1$_{\mathbf{LCD}}$**): i) Solving a sub-problem given $\tau_\mathbf{0}$ as in [24], ii) Our proposed near-optimal solution in Section 3.1 which solves **P1$_{\mathbf{LCD}}$** iteratively, given $\tau_\mathbf{0}$ (Using the approximation in [24]) and given **E** and iii) The same approach in Section 3.1 except for solving **P1$_{\mathbf{LCD}}$**$(\tau_\mathbf{0})$ using the approximation in [31]. This is shown in Figure 2 in which the max-sum throughput is plotted vs. the distance between the AP and ER ($d_{ER-AP}$) for $K = 5$ and $T = 2$. The figure shows the superior performance of our proposed near-optimal solution.

In Figure 3, the max-sum throughput for different values of $K$ at $T = 2$ is plotted vs. the distance between the AP and ER ($d_{ER-AP}$). As noticed, the sum-throughput decreases as the distance increases due to the path loss. As the number of users increases, the sum-throughput increases as well due to existence of more nodes contributing to the sum-throughput. The figure shows that SICD outperforms LCD, which is expected due to the merits of SICD successively decoding interfering signals and canceling them out from the interference to other signals. On the other hand, SICD adds up more computational complexity to the system.

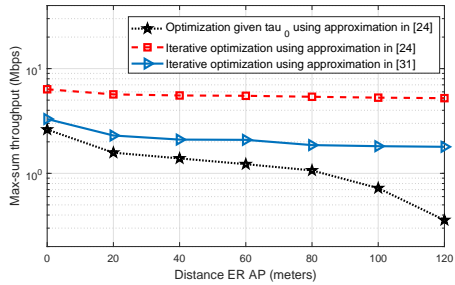
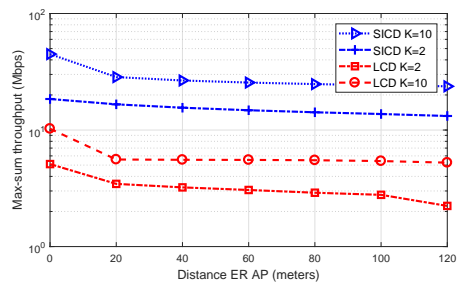

Figure 2: Max-sum at $K = 5$ and $T = 2$

Figure 3: Max-sum at $T = 2$

Moreover, in Figure 4, the max-sum throughput for different values of $T$ at $K = 5$ is plotted vs. the distance between the AP and ER. An increase is observed in the max-sum throughput, as the time horizon, $T$, is extended.

At a $d_{ER-AP}$ of 100 meters, the effect of the decodability threshold, $S_i^{th}$ is studied in Figure 5. **P1$_{\mathbf{LCD}}$** has a solution only for very low values (i.e. low SINR requirements) of the decodability threshold; that's why LCD curves stop early in Figure 5. Moreover, the higher the number of users, the sooner **P1$_{\mathbf{LCD}}$** becomes infeasible, due to elevated levels of interference. On the other hand, **P1$_{\mathbf{SICD}}$** is not affected by the decodability threshold. This is because, with SICD, the average SINR of every user is always greater than the threshold due to interference cancellation.

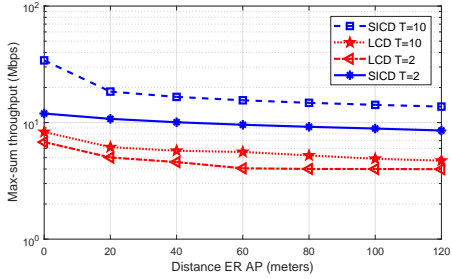
Figure 4: Max-sum at $K = 5$

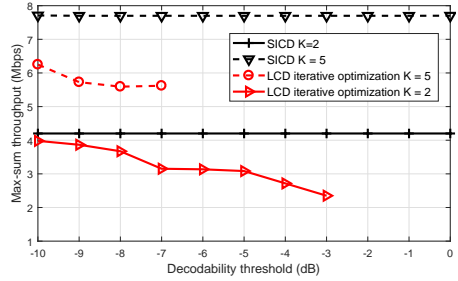
Figure 5: Max-sum at $d_{ER-AP} = 100\ m$

*5.2.2. Max-Min Performance Results*

In an attempt to enhance fairness across users, the objective function to be maximized, under max-min, is the minimum user's throughput. The performance of the presented approach and comparison to the optimum using exhaustive search are presented next.

In Figure 6, the max-min throughput is plotted vs. $d_{ER-AP}$ for both LCD and SICD. The minimum throughput is noticed to decrease as the number of nodes increases. This is because the same resources are distributed among a larger number of users. The performance of SICD is superior to LCD due to higher SINRs as a result of interference cancellation.

To compare the performance of the presented schemes to the global optimum, we resort to exhaustive search since solving for the global optimum is prohibitively complex due to problem non-convexity. Carrying out exhaustive search for a small network scenario ($K = 2$ and $T = 1$), the results are shown in Figure 7.

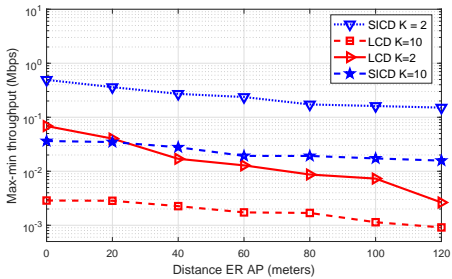
Figure 6: Max-min at $T = 2$

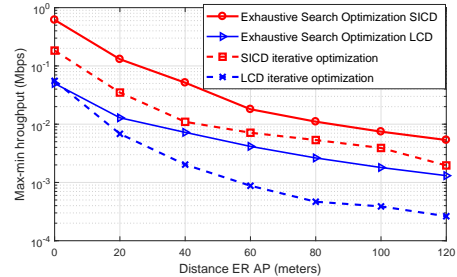
Figure 7: Max-min at $K = 2$ and $T = 1$

*5.2.3. Throughput-Fairness Trade-off*

There is a fundamental trade-off between the max-sum and the max-min problems. A network with max-sum throughput objective will allocate more resources to nodes with high channel gains (near nodes) to add up to the sum-throughput. On the other hand, a network with max-min throughput objective will balance the resource allocation to accommodate users with low channel gains (far nodes) to enhance fairness.

To show this behavior, the sum throughput achieved by the optimal max-min problem is compared to the sum throughput achieved by the max-sum problem. Figure 8 confirms the expected behavior whereby the max-min formulation enhances fairness (as shown in Figures 9 and 10) at the expense of reduced sum-throughput, compared to the max-sum formulation. On the other hand, and in order to complete the picture, the minimum user throughput is compared, under both formulations, in Figure 9 showing that the max-min formulation is superior to max-sum, due to its fairness merits.

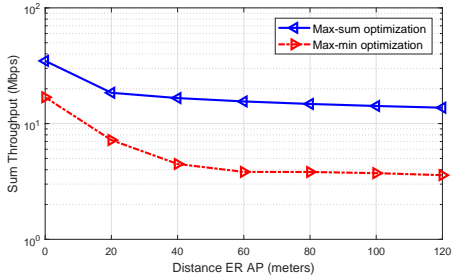 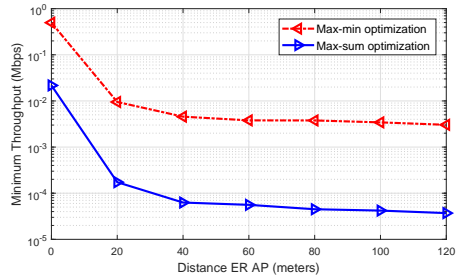

Figure 8: Sum throughput in max-min and max-sum problems

Figure 9: Minimum throughput in max-min and max-sum problems

Moreover, Jain's index (a measure of fairness between nodes) is plotted for the max-sum and max-min optimal policies to show to what extent fairness is accommodated in each formulation.

Jain's index [33] is defined as:

$$j(x_i) = \frac{\left(\sum_{i=1}^{K} x_i\right)^2}{K \cdot \sum_{i=1}^{K} x_i^2}, \qquad (17)$$

where, $x_i$ is the throughput for user $i$, $i$ ranges from 1 to $K$. Jain's index ranges from $\frac{1}{K}$ (worst fairness) to 1 (best fairness), and it is maximum when all users receive the same allocation.

In Figure 10, Jain's index is plotted vs. the number of users. It is noticed how Jain's index in the max-sum problem has actually the worst case value ($\frac{1}{K}$), while it is better for max-min schemes.

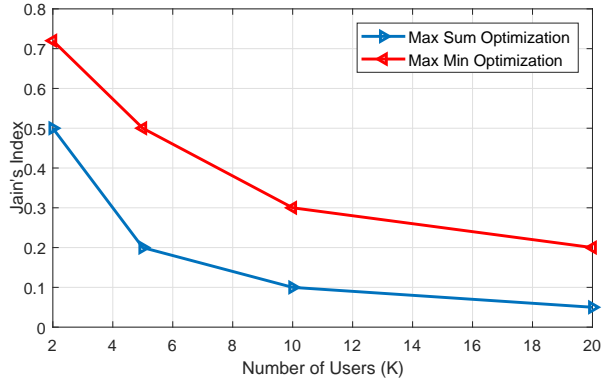

Figure 10: Jain's index (SICD scheme)

Finally, it worth noting that Figures 8, 9 and 10 are plotted for SICD. Similar results can be obtained for LCD.

## 6. Conclusion

In this paper, we investigate the problem of optimal resource (time and power) allocation in WPCNs using non-orthogonal multiple access. Two different optimization problems formulations are considered; in the first one, the sum-throughput (max-sum) of all users is maximized. In the second one, the min-throughput (max-min) of all users is maximized. Under these two formulations, two NOMA decoding schemes are studied, namely, low complexity decoding (LCD) and successive interference cancellation decoding (SICD). Due to the non-convex nature of the max-sum and max-min optimization problems, we propose an approximate solution approach, in which the non-convex optimization problem is approximated by a convex optimization problem, which satisfies all the constraints of the original problem. The approximate convex optimization problem can be then solved iteratively. The results show a trade-off between maximizing the sum throughout and achieving fairness via maximizing the minimum throughput.